\documentclass[12pt]{article}

\usepackage{graphicx}
\usepackage{amsmath,amsthm,amssymb,color}
\usepackage[noadjust,sort]{cite}
\usepackage[mathscr]{euscript}
 \usepackage[normalem]{ulem} 
\usepackage{subcaption}
\usepackage{xcolor}
\usepackage{float}
\usepackage{bbm}

\usepackage[font=small,labelfont=bf]{caption}
\usepackage{framed}
 \usepackage[boxruled]{algorithm2e}

\normalsize

{\endMakeFramed}

\theoremstyle{plain}
\newtheorem{Theorem}{Theorem}[section] %

\theoremstyle{definition}
\newtheorem{Remark}[Theorem]{Remark}

\def\naive{f_{\operatorname{naive}}}
\def\calF{\mathcal{F}}

\def\calL{\mathcal{L}}
\def\calH{\mathcal{H}}

\def\calR{\mathcal{R}}
\def\Err{\mathfrak{E}}

\theoremstyle{definition}

\newtheorem{Problem}{Problem}[section]

{\par\noindent{\it Proof of}} 
{\hfill$\vspace{5mm}\scriptstyle\blacksquare$} 

\numberwithin{equation}{section} 

\def\supp{\operatorname{supp}}

\def\R{\mathcal{R}}

\def\Err{\mathfrak{E}}

\def\Reals{{\mathbb{R}}}

\def\Naturals{{\mathbb{N}}}

\def\st{\,:\,}

\def\dfrac#1#2{\lower0.15ex\hbox{\large$\textstyle\frac{#1}{#2}$}}

\def\leq{\leqslant}
\def\geq{\geqslant}

\begin{document}

\setcounter{page}{1}

\title{
	PSWF-Radon approach to 
reconstruction from band-limited Hankel transform
%
}
\date{}

\author{ 
\and
Fedor  Goncharov\\
\small Université Paris-Saclay\\[-0.8ex] 
\small CEA, List, F-91120\\[-0.8ex]
\small Palaiseau, France
\\[-0.8ex]
\small\texttt{fedor.goncharov@cea.fr}
\and
Mikhail Isaev
\\\small
School of Mathematics 
and Statistics
 \\[-0.8ex]
\small UNSW Sydney\\[-0.8ex]
\small Sydney, NSW, Australia\\[-0.8ex]
\small\tt  isaev.m.i@gmail.com  
\and
Roman G. Novikov
\\
\small CMAP, CNRS, Ecole polytechnique\\[-0.8ex]
\small Institut Polytechnique de Paris\\[-0.8ex]
\small Palaiseau, France\\
\small IEPT RAS, Moscow, Russia\\
\small\texttt{novikov@cmap.polytechnique.fr}
\and
Rodion Zaytsev \\
\small
Faculty of Mathematics, HSE University
\\
\small
Igor Krichever Center for Advanced Studies 
\\[-0.8ex]
\small 
Skolkovo Institute of Science and Technology
\\
\small Moscow, Russia
\\
\small\texttt{rvzaytsev@edu.hse.ru}
}

\maketitle

\begin{abstract}
	We give new formulas for reconstructions from band-limited Hankel transform of  integer or half-integer order. 
	Our formulas rely on the PSWF-Radon approach to super-resolution in multidimensional Fourier analysis.  This approach consists of combining the theory of classical one-dimensional prolate spheroidal wave functions with the Radon transform theory.  We also use the relation between Fourier and Hankel transforms and  
	Cormack-type inversion of the Radon transform. 
Finally, we investigate numerically 
 the capabilities of our approach to super-resolution for band-limited Hankel inversion
in relation to varying levels of noise.
 \\
	
\end{abstract}

\section{Introduction}\label{S:intro}
The Hankel transform  of order $\nu$ is  formally defined by
\begin{equation}\label{def:Hankel}
 \calH_{\nu}[f](t) 
:= \int_0^{\infty} f(s) J_{\nu}(ts) \sqrt{ts}\, d s, \qquad t \in \Reals_+,
\end{equation}
where $J_{\nu}$ is the Bessel function of the first kind and 
 $\Reals_+:=\{\rho \in \Reals \st \rho\geq 0\}$.
If $\nu\geq -\frac{1}{2}$ then $\calH_\nu$ is an invertable operator on $\calL^2(\Reals_+)$ and
  $\calH_{\nu}^{-1} = \calH_{\nu}$, that is, 
 the operator defined by \eqref{def:Hankel} coincides with  its own inverse. 
%
 

We consider the following inverse problem.
\begin{Problem}\label{P1}
	Let  $\sigma, r>0$ and $\nu \geq -\frac 12$ be given.  Find   $f \in \calL^2(\Reals_+)$ from  
		 $h=\calH_{\nu}[f] $  given on $[0,r]$ (possibly with some noise), 
		 under a priori assumption that $\supp f \subseteq  [0,\sigma]$.
	\end{Problem}

Problem \ref{P1} is uniquely solvable for the   noiseless case in view of the analyticity  of $\calH_{\nu}[f]$  on $(0,+\infty)$  under our assumptions.  On the other hand, this problem is severely ill-posed in the sense of Hadamard  in view of the super-exponential decay of the eigenvalues of the band-limited Hankel transform; see \cite[Lemma 1]{KM2009}.
For introduction into the theory of ill-posed problems; see, for example, 
 \cite{HR2021,TA1977}.


In the present work, we consider Problem \ref{P1} with integer or  half-integer order $\nu \geq 0$ in connection with  the   problem of reconstructing 
 unknown function   supported in  the ball $B_\sigma\subset \Reals^d$, $d\geq 2$, from its Fourier transform given on the ball $B_r$ under a priori assumptions 
 of the cylindrical or spherical symmetry, and also for more general cases of variable separation in polar coordinates.
  Here, for $\rho \in \Reals_+$,  we let 
 \[
 	   B_\rho :=\{ q \in \Reals^d \st |q|\leq \rho\}.
 \]
Applications of Problem \ref{P1} include   multidimensional monochromatic inverse scattering in the Born approximation and inverse source problems; see
 \cite{INS2022, IN2020+, Meng2023, AM, NS2024}
and references therein.  See also  \cite{ KB2006,  KM2009}  for other important applications.  
For information on the problem of band-limited Fourier inversion 
without  a priori assumptions of spherical-type symmetry, including various results and approaches; see, for example, \cite{AMS2009, AM, BM2009,  CF2014, Gerchberg1974, IN2020, IN2020+,INnotePSWF, INS2022,INS_Proceedings, LRC1987, Meng2023, Papoulis1975}.

A possible  approach  for solving Problem  \ref{P1}   is based on the following approximation
\begin{equation}\label{eq:naive}
	f \approx f_{\rm naive} :=  \calH_{\nu}^{-1} \left[h^{\rm ext} \right]  \text{ on } [0,\sigma],
\end{equation}
where  
\[
	h^{\rm ext}(t):= 
	\begin{cases}
			h (t), &  \text{ for } t \in [0,r],\\
			0, & \text{otherwise.}
		\end{cases}
\]
Formula \eqref{eq:naive} leads to a stable and accurate reconstruction for  sufficiently   large $r$.   However, similarly to the reconstruction from the Fourier  transform given on the ball $B_r$,  this naive approach has   the   \emph{diffraction limit}: small details  (especially less than  $\pi /r$) are blurred.  The term \emph{super-resolution} refers to the techniques that allow reconstruction beyond this diffraction limit.

%
%
%
%
%
%

 In the present work we give new  theoretical and numerical results  for  Problem \ref{P1} with integer or half-integer $\nu\geq 0$ following   the PSWF-Radon approach of \cite{INnotePSWF, INS2022, INS_Proceedings} to Fourier analysis in dimension  $d=2,3$. 
%
%
In connection with the PSWF theory, we consider  the  integral operator $\calF_c$ on $\calL^2([-1,1])$, defined by
\begin{equation}\label{def:Fc}
	\calF_c[f] (x) := \int_{-1}^1 e^{i c xy} f(y)dy,
\end{equation}
    where  $c>0$ is the \emph{bandwidth parameter}.   
   The eigenfunctions 
   $(\psi_{j,c})_{j \in \Naturals}$    of $\calF_c$ 
   are   \emph{prolate spheroidal wave
   	functions} (PSWFs), where  $\Naturals:=\{0,1\ldots\}$.
   Recall that  the operator $\calF_c$ and its inverse admit the singular value decompositions:
      \begin{align}
      	\calF_c [f] (x) &= \sum_{j \in \Naturals} \mu_{j,c}\psi_{j,c}(x) \int_{-1}^1 \psi_{j,c} (y) f(y) dy \label{Fc-dec}, \\
   	\calF_{c}^{-1} [g](y) &=    \sum_{j \in \Naturals} \dfrac{1}{\mu_{j,c}}\psi_{j,c}(y) \int_{-1}^1 \psi_{j,c} (x) g(x)dx.
   	\label{f:inverse}
   \end{align}
 In addition, we can assume that   
   $0<|\mu_{j+1,c}| < |\mu_{j,c}|$ for all $j \in \Naturals$. For more details on the PSWF theory,  see \cite{AGD2014,BK2014, BK2017, INnotePSWF, KRD2021, Wang2010, XRY2001, SP1961}.

First, we consider the case of integer $\nu$.

\begin{Theorem}\label{T:int}
	Let $\nu \in \Naturals$, $r,\sigma >0$, and $c:= r\sigma$.	Let  $f\in \calL^2(\Reals_+)$ be supported in $[0,\sigma]$.
Then,   transform $h:= \calH_\nu[f]$ on $[0,r]$  uniquely determines $f$ by the formula 
\begin{equation} \label{f:int}
f(s)
= 
-   \frac{2 i^\nu  }{\sigma}   \sqrt{s}\frac{d}{   ds} 
\int_s^{\sigma} \frac{ s \, T_{\nu}\left(t/s\right)}{  t(t^2 - s^2)^{\frac{1}{2}}}   \calF_{c}^{-1} [h_{r,\nu}](t/\sigma)     dt, \qquad s \in (0,\sigma],
\end{equation}
where $T_{\nu}$ is the  Chebyshev polynomial of the first kind of order $\nu$ and 
	\begin{equation}\label{g:int}
	h_{r,\nu}(x) := \begin{cases} \frac{1}{\sqrt{rx}}h (rx), &  \text{if $x \in [0,1] $,}\\
		(-1)^{\nu}  \frac{1}{\sqrt{-rx}} h(-rx), & \text{if $x \in [-1,0) $.}
	\end{cases}
	\end{equation}
	\end{Theorem}

Second, we consider the case of half-integer $\nu$.

\begin{Theorem}\label{T:half}
	Let $n := \nu-\frac12 \in \Naturals$, $r,\sigma >0$, and $c := r\sigma$. Let  	$f\in \calL^2(\Reals_+)$ be supported in $[0,\sigma]$.
Then,     $h:=\calH_\nu[f]$ on $[0,r]$  uniquely determines $f$ by the formula 
	\begin{equation}\label{f:half}
f(s)
= 
  \frac{\sqrt{2 \pi}\,i^n }{\sigma}   s \frac{d^2}{ds^2} 
\int_s^{\sigma} \frac{s P_{n}\left(t/s\right)}{t^2}   \calF_{c}^{-1} [h_{r,\nu}](t/\sigma)   dt, \qquad s \in (0,\sigma],
\end{equation}
where $P_{n}$ is  the  Legendre polynomial of order $n$ and 
	\begin{equation}\label{g:half}
	h_{r,\nu}(x) := \begin{cases} \frac{1}{rx}h (rx), & \text{if $x \in [0,1] $,}\\
		(-1)^{n}  \frac{1}{rx} h(-rx), &\text{if $x \in [-1,0) $.}
	\end{cases}
	\end{equation}
	\end{Theorem}

Note that  the function $h_{r,\nu}$ considered in \eqref{g:int}  and \eqref{g:half} 
has no singularity at $0$ under our assumptions, in view of 
 the definition  of the Hankel transform in \eqref{def:Hankel}.

 We prove Theorem \ref{T:int} and Theorem \ref{T:half}   in Section \ref{S:proofs} using
the PSWF-Radon approach for reconstruction from the band-limited Fourier transform  developed in  \cite{INnotePSWF, INS2022, INS_Proceedings},
 classical formulas relating the Fourier   and  Hankel transforms,  and  Cormack-type formulas    for inversion of the Radon transform.
  It is known that Cormack's formulas have good numerical properties for small $\nu$ but instability increases exponentially as $\nu$ grows; see, for example, \cite[Section II.2]{Naterrer2001}. 
  A more stable approach relies on  the
  \emph{filtered back projection} (FBP) algorithm which  can be derived from the classical Radon inversion formula; 
  see   \cite{Puro2001,Naterrer2001}   for more details
  in the case  of  $d=2$ (and the same  ideas would work for $d=3$).
  Therefore, we employ FBP-based algorithms for our numerical implementation of  band-limited Hankel transform inversion; see Section \ref{S:numerics} for detailed presentation.

\begin{figure}[h!]
\centering

 \begin{subfigure}{0.32\textwidth}
 \includegraphics[height=5cm]{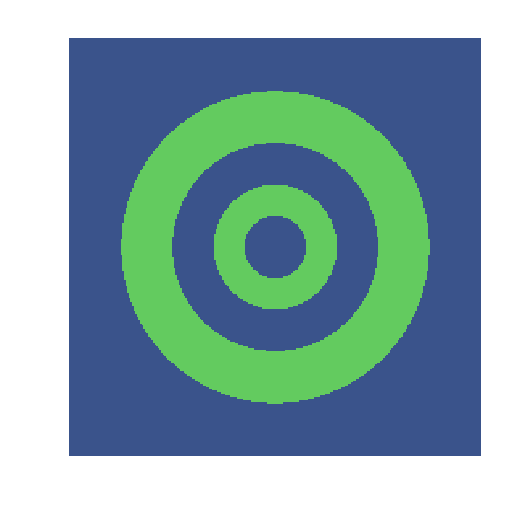}
 \caption{
 }
\end{subfigure}
 \begin{subfigure}{0.32\textwidth}
 \includegraphics[height=5cm]{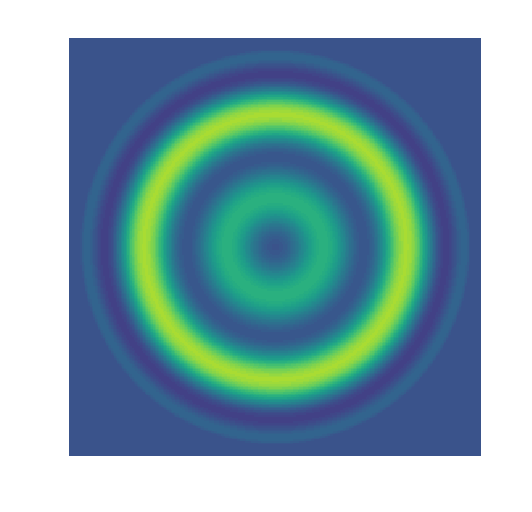}
 \caption{
 }
\end{subfigure}
 \begin{subfigure}{0.32\textwidth}
 \includegraphics[height=5cm]{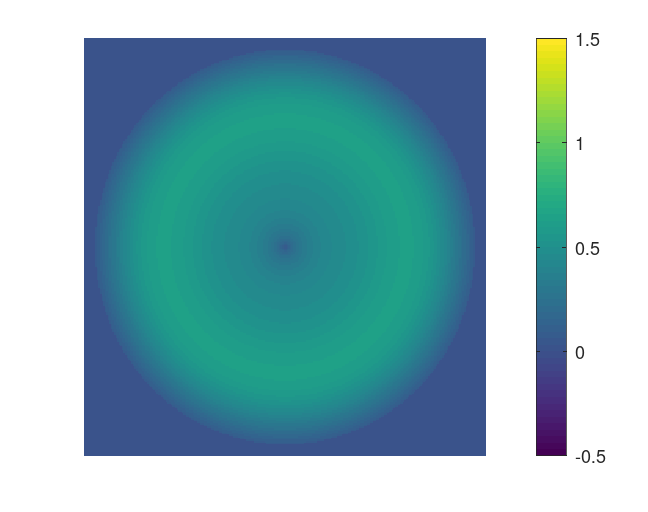}
 \caption{
 }
\end{subfigure}
\captionsetup{singlelinecheck=off}
 \caption[foo bar]{  \label{Fig2D}
  An example of reconstruction for Problem \ref{P1} with $\nu=0$, $\sigma = 1$,  $r=10$: 
  (a) preimage, (b) PSWF-Radon reconstruction, (c)  reconstructio via \eqref{eq:naive}. 
 }
\end{figure}

\begin{figure}[h!]
\centering

 \begin{subfigure}{0.32\textwidth}
 \includegraphics[height=5cm]{phantom.png}
 \caption{
 }
\end{subfigure}
 \begin{subfigure}{0.32\textwidth}
 \includegraphics[height=5cm]{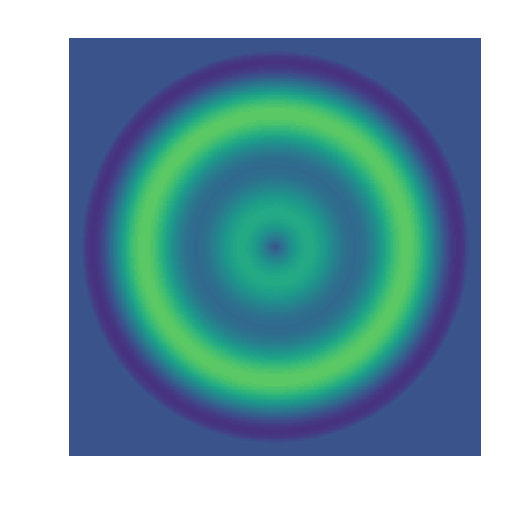}
 \caption{
 }
\end{subfigure}
 \begin{subfigure}{0.32\textwidth}
 \includegraphics[height=5cm]{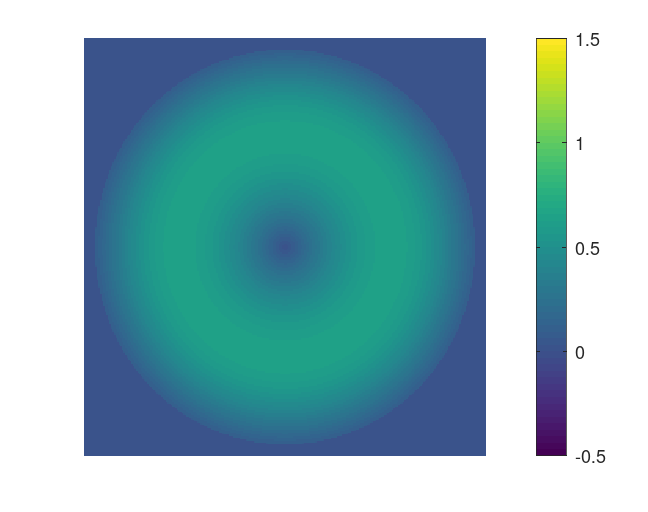}
 \caption{
 }
\end{subfigure}
\captionsetup{singlelinecheck=off}
 \caption[foo bar]{ \label{Fig3D}
  An example of reconstruction for Problem \ref{P1} with $\nu=0.5$,
  $\sigma = 1$,  $r=10$:
  (a) preimage, (b) PSWF-Radon reconstruction, (c)~reconstruction via \eqref{eq:naive}. 
 }
\end{figure}

  Figure \ref{Fig2D} and   Figure \ref{Fig3D} illustrate our numerical reconstruction for Problem \ref{P1}, where  $\nu= 0$ and $\nu=0.5$, respectively, and $\sigma = 1$,  $r=10$.  Here, we use two-dimensional visialisation of our results  in the form $v(x) = f(|x|)$, where  $x \in \Reals^d$ for $d=2$ and for $d=3$ (displaying central cross section), for the preimage,  its naive reconstruction based on \eqref{eq:naive}, and its PSWF-Radon reconstruction with an appropriate regularisation. 
  In particular, this also illustrates reconstruction of $v$   from its  band-limited Fourier transform under additional a priori assumption of spherical symmetry (in view of classical relation between the Fourier and Hankel transforms).
 The preimage considered in Figures \ref{Fig2D} and~\ref{Fig3D} is the characteristic function of   $[0.15,0.3]\cup [0.5,0.75]$. Note that the gaps and the internal ring width  are essentially smaller that   the diffraction limit $\pi/10\approx 0.31$; therefore,  the standard naive approach via \eqref{eq:naive} is incapable to distinguish these details. In contrast, the PSWF-Radon approach succeeds to recover preimage for both $\nu =0$ and $\nu = 0.5$. 
Furthermore, for this preimage, 
our super-resolution reconstruction 
remain  adequate   even in the presence of moderate noise of $10\%$ for $\nu=0$ and $5\%$ for $\nu=0.5$;
 see  Section \ref{S:numerics} for more detailed investigation on the stability of the PSWF-Radon approach to the band-limited Hankel transform inversion.

  

In Section 4, we provide additional discussion on the PSWF-Radon approach,   comparing it with other possible approaches to Problem \ref{P1} and outlining future directions for development.


\section{Proofs of the theorems}\label{S:proofs}

We consider the classical Fourier transform $\calF$ defined by
\begin{equation}\label{def:Fourier}
    \calF[v](p) := \frac{1}{(2\pi)^d}
    \int_{\Reals^d} e^{ipq} v(q) dq, \qquad p\in \Reals^d.
    \end{equation}
First, we
recall the formulas  from  \cite{INnotePSWF}  for reconstruction $v$ from $\calF[v]$  limited to the ball $B_r$ for a compactly supported $v$.
Then, we prove Theorem \ref{T:int} and Theorem \ref{T:half}.

 \begin{Theorem}[Isaev, Novikov \cite{INnotePSWF}]\label{T:Fourier}Let $d\geq 1$,  $r,\sigma>0$ and $c = r\sigma$.
  Let $v\in \calL^2(\Reals^d)$ and $\supp v \subset  B_{\sigma}$. 
  Then,   its Fourier transform $\calF[v] $ restricted to $B_r$ determines $v$ via the following formulas:
  \begin{align*}
  v(q) &=     \R^{-1} [w_{r, \sigma}]( \sigma^{-1}q), \qquad  q\in \Reals^d,\\
  w_{r,\sigma}(y,\theta) &=
 	 \begin{cases}
 \calF_c^{-1}[g_{r,\theta}](y), &\text{if } y\in [-1,1]\\
 0, & \text{otherwise},  
 \end{cases} 
 \\
 g_{r,\theta} (x)&=      \left(\dfrac{2\pi}{\sigma}\right)^d \calF[v] (r x\theta), \qquad  x\in [-1,1],\ \theta \in \mathbb{S}^{d-1},
  \end{align*}
 where     $\calF_c^{-1}$ is     defined by \eqref{f:inverse}
 and   $ \R^{-1}$ is the inverse Radon transform.
 \end{Theorem}

In the above theorem, $\calR$ is the classical Radon transform defined by 
\[
    \mathcal{R}_{\theta}[v](y)= w(y,\theta):= \int_{q\in \Reals^d \st q\theta = y} v(q) dq, \qquad y\in \Reals,\, \theta \in \mathbb{S}^{d-1}. 
\]
There are multiple ways to compute  $\calR^{-1}$; see, for example, \cite{Naterrer2001}.  
We use  Cormack-type formulas assuming that $v$ is compactly supported. For    $d=2$ and $q\neq 0$,  we have 
\begin{equation}\label{Cormack2D}
\begin{aligned}
   v(q) &= \calR^{-1}[w] = \sum_{n=-\infty}^{+\infty} v_n(|q|) e^{i n \phi}, \qquad \frac{q}{|q|}= (\cos \phi, \sin \phi)\\
   v_n(s) &= -\frac{1}{\pi}\frac{d}{ds} \int_{s}^{+\infty} \frac{s T_{|n|}(t/s) w_n(t) }{ t \sqrt{t^2-s^2} } dt, \qquad s>0,
   \end{aligned}
 \end{equation} 
   where  $T_{|n|}$ is the Chebyshev polynomial of the first kind and  $(w_n)_{n\in \mathbb{Z}}$ are defined by
   \[
   w(t,\theta) = \sum_{n=-\infty}^{+\infty}   w_n(t) e^{in\psi}, \qquad t>0, \ \theta = (\cos\psi, \sin \psi).
\]
For   $d=3$ and $q\neq 0$, we have 
\begin{equation}\label{Cormack3D}
       \begin{aligned}
        v(q)&= \calR^{-1}[w] = \sum_{n=0}^{+\infty} \sum_{k = -n}^{n}  v_{n,k}(|q|) \,Y_{n,k}\left(\dfrac{q}{|q|}\right),
\\
 v_{n,k}(s) &= \frac{1}{2\pi}  \cdot \frac{d^2}{ds^2} \int_{s}^{+\infty} \frac{s P_{n}(t/s) }{t^2}  w_{n,k}(t) dt, \qquad s>0,
 \end{aligned}
\end{equation}
where $Y_n^k:\mathbb{S}^2 \rightarrow \Reals$ is the spherical harmonic for each $n \in \mathbb{N}$ and $|k|\leq n$,   $P_n$ is  the  Legendre polynomial, and $ (w_{n}^k)_{n \in \mathbb{N}, |k|\leq n}$ are defined by
\[
    w(t,\theta) = \sum_{n=0}^{+\infty} \sum_{k = -n}^{n}  w_{n,k}(t) Y_{n,k}(\theta), \qquad t>0, \ \theta \in \mathbb{S}^2.
\]

\begin{Remark}
     Note that formula \eqref{Cormack3D} is usually stated in an alternative form with
     \[
      v_{n,k}(s)  =  \frac{1}{2\pi s}\int_{s}^{+\infty}   P_{n}(t/s)   w_{n,k}''(t) dt, \qquad s>0.
     \]
     where $w_{n,k}''$ is the second derivative of $w_{n,k}$;
     see, for example, \cite[Theorem 2.3]{Naterrer2001}. 
     If $v$ is sufficiently smooth and compactly supported (implying that  $w_{n,k}$ is such as well), one can show the equivalence by  using the following formula twice:  
     \[
        \int_s^{+\infty} F(t/s) G'(t) dt  = 
        s \frac{d}{ds} \left(\frac{1} {s}\int_{s}^{+\infty}  \tilde{F}(t/s)  G(t) dt\right), \qquad \tilde{F}(t/s)  := \frac{t}{s}F(t/s),
     \]
    where $F,G $ are smooth on $(0,+\infty)$ and $G$ is  compactly supported. 
    This also extends to non-smooth $v$, but we prefer the form of \eqref{Cormack3D} as it does not require  consideration of generalised derivatives.  
     \end{Remark}

\subsection{Proof of Theorem \ref{T:int}}
For non-zero $q \in \Reals^2$, we introduce the notation $\varphi_q \in [0,2\pi)$  defined by
\[
    \frac{q}{|q|}= (\cos \varphi_q, \sin \varphi_q).
\]
Recall that if    $v \in \calL^2(\Reals^2)$  is expanded in the   series
\begin{equation*} 
    v(q) = \sum_{n=-\infty}^{+\infty} v_n(|q|) e^{in\varphi_q},
\end{equation*}
then its Fourier transform \eqref{def:Fourier} satisfies 
\begin{equation}\label{relation_2D}
    \calF[v](p) = \frac{1}{2\pi \sqrt{|p|}}\sum_{n=-\infty}^{+\infty} i^{n}  \calH_n[f_n](|p|) e^{in \varphi_p},  
    \end{equation}
    where $p \in \Reals^2$ is non-zero and
    \[
    f_n(s) := v_n(s)\cdot \sqrt{s}, \qquad s\in \Reals_+.
\]

To prove Theorem \ref{T:int}, we 
apply Theorem \ref{T:Fourier} to the function $v$
defined by
 \[
  v(q) := \frac{f(|q|)}{\sqrt{|q|}} e^{i\nu\varphi_q}.
\]
Since, by our assumptions, $f\in \calL^2(\Reals_+)$,  we get that $v\in \calL^2(\Reals^2)$.
Using \eqref{relation_2D}, we get that 
\[
     \calF[v] (p) = 
   \frac{i^\nu}{2\pi \sqrt{|p|}}\calH_{\nu}[f] (|p|)  e^{i\nu \varphi_p}.
\]
Taking $p = rx \theta$ and recalling $h_{r,\nu}$ 
defined in \eqref{g:int}, we find  that 
\begin{equation}\label{g_2D}
    g_{r,\theta}(x):= \left(\frac{2\pi}{\sigma}\right)^2 \calF[v] (rx  \theta) = \frac{ 2\pi i^{\nu} }{\sigma^2 }h_{r,\nu}(x) e^{i \nu \varphi}, \qquad x\in [-1,1]. 
\end{equation}
where $\theta = (\cos \varphi, \sin \varphi)$. Next, we observe that 
\[
    w_{r,\sigma}(y,\theta) :=  \calF_c^{-1}[g_{r,\theta}](y) = \frac{ 2\pi i^{\nu} e^{i \nu \varphi} }{\sigma^2  } \calF_c^{-1}[ h_{r,\nu}](y) , \qquad y\in [-1,1].
\]
Applying  Theorem \ref{T:Fourier} and Cormack's formula \eqref{Cormack2D}, we find that 
\begin{align*}
    v( \sigma x \theta) &= \calR^{-1}[w_{r,\sigma}] (x)
    \\
    &=
    -
    \frac{e^{i \nu \varphi}}{\pi}\frac{d}{dx} \int_{x}^{+\infty} \frac{x T_{\nu}(y/x) w_\nu(y) }{ y \sqrt{y^2-x^2} } dy, \qquad x\in (0,1],
\end{align*}
where 
\[
    w_\nu(y):= 
    \begin{cases}
    \dfrac{ 2\pi i^{\nu} }{\sigma^2  } \calF_c^{-1}[ h_{r,\nu}](y), &\text{if $y\in (0,1]$,}
    \\
    0, &\text{if $y>1$.}
    \end{cases}
\]
Substituting  $s = \sigma x$, $t= \sigma y$, and recalling the definition of $v$, we derive formula \eqref{f:int}. This completes the proof of Theorem \ref{T:int}.

\subsection{Proof of Theorem \ref{T:half}}
The proof follows  the same steps as in the previous section.
First, we recall that if    $v \in \calL^2(\Reals^3)$  is expanded in the   series
\begin{equation*} 
    v(q) = \sum_{n=0}^{+\infty} \sum_{k = -n}^{n}  v_{n,k}(|q|) \,Y_{n,k}\left(\dfrac{q}{|q|}\right),
\end{equation*}
then its Fourier transform \eqref{def:Fourier} satisfies 
\begin{equation}\label{relation_3D}
    \calF[v](p) = \frac{1}{(2\pi)^{3/2} |p|}\sum_{n=0 }^{+\infty} \sum_{k = -n}^{n}  i^{n}  \calH_{n+\frac12}[f_{n,k}](|p|) Y_{n,k}\left(\dfrac{p}{|p|}\right),
    \end{equation}
    where $p\in \Reals^3$ is non-zero and 
    \[
    f_{n,k}(s) := v_{n,k}(s)\cdot s, \qquad s\in \Reals_+.
\]

 To prove Theorem \ref{T:half}, we 
apply Theorem \ref{T:Fourier} to the function $v$
defined by
 \[
  v(q) := \frac{f(|q|)}{ |q| } Y_{n,0}\left(\dfrac{q}{|q|}\right).
\]
Since, by our assumptions, $f\in \calL^2(\Reals_+)$,  we get that $v\in \calL^2(\Reals^3)$. Using \eqref{relation_3D}, we get that 
\[
     \calF[v] (p) = 
   \frac{i^n}{ (2\pi)^{3/2}  |p| }\calH_{\nu}[f] (|p|)  Y_{n,0}\left(\dfrac{p}{|p|}\right).
\]
Taking $p = rx \theta$  and recalling $h_{r,\nu}$ 
defined in \eqref{g:half}, we find  that 
\begin{align*}
    g_{r,\theta}(x):&= \left(\frac{2\pi}{\sigma}\right)^3 \calF[v] (rx  \theta) \nonumber \\ &= \frac{ (2\pi)^{3/2} i^{n} }{\sigma^3 }h_{r,\nu}(|x|) Y_{n,0}\left(\operatorname{sgn}(x)\theta\right), \quad  x\in [-1,1],\quad  \theta\in \mathbb{S}^2. 
\end{align*}
Recall that the spherical harmonic $Y_{n,k}$ with $k=0$ can be expressed in terms of the   Legendre polynomial  $P_n$ and, in particular,
$   
    Y_{n,0}\left(-\theta\right) 
    = (-1)^n Y _{n,0}\left(\theta\right).
$ 
Therefore, we get
\begin{equation}\label{g_3D}
  g_{r,\theta}(x) = \frac{ (2\pi)^{3/2} i^{n} }{\sigma^3 }h_{r,\nu}(x) Y_{n,0}\left(\theta\right)
\end{equation}
for all non-zero $x\in [-1,1]$.  Next, we observe that 
\[
    w_{r,\sigma}(y,\theta) :=  \calF_c^{-1}[g_{r,\theta}](y) = \frac{ (2\pi)^{3/2} i^{n}   }{\sigma^3  } 
    \calF_c^{-1}[ h_{r,\nu}](y)  Y_{n,0}\left(\theta\right), \qquad y\in [-1,1].
\]
Applying  Theorem \ref{T:Fourier} and the Cormack-type formula of \eqref{Cormack3D}, we find that 
\begin{align*}
    v( \sigma x \theta) &= \calR^{-1}[w_{r,\sigma}] (x)
    \\
    &=
    \frac{Y_{n,0}\left(\theta\right)}{2\pi}\frac{d^2}{dx^2} \int_{x}^{+\infty} \frac{x P_{n}(y/x) w_n(y) }{  y^2 } dy, \qquad x \in (0,1],
\end{align*}
where 
\[
    w_n(y):= 
    \begin{cases}
    \dfrac{ (2\pi)^{3/2} i^{n} }{\sigma^3  } \calF_c^{-1}[ h_{r,\nu}](y), &\text{if $y\in (0,1]$,}
    \\
    0, &\text{if $y>1$.}
    \end{cases}
\]
Substituting  $s = \sigma x$, $t= \sigma y$, and recalling the definition of $v$, we derive formula \eqref{f:half}. This completes the proof of Theorem \ref{T:half}.

 
\section{Reconstructions via the FBP algorithms}\label{S:numerics}
Numerical reconstruction for  Problem \ref{P1} with integer or half-integer $\nu\geq 0$ can be implemented as follows.  As  discussed in Section \ref{S:intro}, we rely on \emph{filtered back projection} (FBP) algorithms
for   inverting the Radon transform,  unlike Theorem~\ref{T:int} and Theorem \ref{T:half} which rely on Cormack-type inversion.

\begin{enumerate}
\item[Step 1.] 
Given the Hankel data $h$ on $[0,r]$, 
set $h_{r,\nu}$  according formulas \eqref{g:int} and \eqref{g:half} as in Theorem~\ref{T:int} and Theorem \ref{T:half}. 
 Define   $g_{r,\theta}$ according to formulas \eqref{g_2D} and \eqref{g_3D}:
\[
    g_{r,\theta}(x) = 
    \begin{cases}
    \dfrac{ 2\pi i^{\nu} }{\sigma^2 }h_{r,\nu}(x) e^{i \nu \varphi}, &\text{if $\nu$ is integer,} \\
     \dfrac{ (2\pi)^{3/2} i^{n} }{\sigma^3 }h_{r,\nu}(x) Y_{n,0}\left(\theta\right), &\text{if $\nu$ is half-integer.} 
    \end{cases}
\]
Recall that $x\in [-1,1]$ and  $n = \nu-\frac 12$ for half-integer $\nu$.

\item[Step 2.]
Reconstruct $v$ from $g_{r,\theta}$ according to the formulas of Theorem \ref{T:Fourier}.
For inverting  $ \calF_c^{-1}$, we use its regularised version 
$\calF_{m,c}^{-1}$ 
as in \cite{INnotePSWF, INS2022, INS_Proceedings},
which we recall below.
For  inverting the Radon transform, we use  the   FBP algorithms; see, for example,   \cite{INS2022, INS_Proceedings, Goncharov2021}.

\item[Step 3.] Compute the result $f$ on $[0,\sigma]$  as follows:
\[
        f(s) = 
     \begin{cases}
     \displaystyle
     \frac{\sqrt{s}}{2\pi}\int_{-\pi}^{\pi} v(s\cos\phi, s\sin \phi) e^{-i\nu \phi} d\phi, &\text{if $\nu$ is integer,}\\
      \displaystyle s \int_{\mathbb{S}^2} v(s\theta) Y_{n,0}(\theta)d\theta, &\text{if $\nu$ is half-integer.} 
     \end{cases}
\]

\end{enumerate}






     
             




 At Step 2,
 we approximate the operator
 $ \calF_c^{-1}$ in \eqref{f:inverse} by the finite rank operator  $\calF_{m,c}^{-1}$ defined by
 \begin{equation}\label{def:Fnc}   
 	\calF_{m,c}^{-1} [g] (y) :=  \sum_{j=0}^m \dfrac{1}{\mu_{j,c}}\psi_{j,c}(y) \int_{-1}^1 \psi_{j,c} (x) g(x)dx. 
 \end{equation}
 The parameter $m$ in \eqref{def:Fnc} 
 is the regularisation parameter; see   
\cite[Section 2]{INS2022} for a discussion 
on its optimal choice.

Let 
$\tilde \calH_{m,\nu}^{-1} \left[h \right]$ denote the result of our reconstruction  described above from the Hankel data $h$ with regularisation parameter $m$.
We use the following measure of the quality of reconstruction $f_{\rm rec}$ of the preimage $f$:
\begin{equation}\label{def:Err}
		\Err(f_{\rm rec}, h)   :=   \frac{\|\calH_{\nu}[f_{\rm rec}]-h\|_{\calL^2([0,r])} }{\|h\|_{\calL^2([0,r])}}.
\end{equation}
In the numerical examples of the present work 
the choice of $m$ based on \emph{the  residual minimisation  principle} (that is,
taking $m$ minimising $\Err(\tilde \calH_{m,\nu}^{-1} \left[h \right], h)$) is highly satisfactory: it  gives a stable reconstruction even from   noisy data (without blow-ups in the configuration space) and leads to super-resolution for moderate levels of noise.   

We present our examples  in Section \ref{S:two-step} and Section \ref{S:harmonics}, where 
we always  have $\sigma=1$ and $c=r=10$.
Our numerical implementations rely on the values of the Hankel data
 $h$  on the standard uniform grid with $N$ points, where we take $N=256$.
In the examples with noise, the values of  $h$
at these points are altered by centered Gaussian white random noise independently of each other. 
 All figures   use the same legend: 
 \begin{itemize}
 \item \emph{dotted lines} represent the preimage $f$ and its Hankel data $h$ (possibly with noise);
 
 \item \emph{bold lines} represent the PSWF-Radon reconstruction $\tilde f = \tilde{f}_m = \tilde \calH_{m,\nu}^{-1} \left[h \right]$  and its   Hankel transform
   with $m$ chosen according to the  residual minimisation  principle (unless specified otherwise); 
 
 \item \emph{dashed lines} represent    $\naive$  defined in \eqref{eq:naive} and its  Hankel transform $\calH_{\nu}[\naive]$.
\end{itemize}
We focus on the cases when $\nu=0$ and $\nu=0.5$ as our reconstruction behave similarly for higher integer and half-integer orders.  As mentioned in Section \ref{S:intro},  this also illustrates reconstructions   from  band-limited Fourier transform under additional a priori assumption of spherical symmetry in dimensions $d=2$ and $d=3$; see Figure \ref{Fig2D} and Figure \ref{Fig3D}.

\subsection{Two step function}\label{S:two-step} 
We start from a more detailed analysis of the example 
considered in Figure~\ref{Fig2D} and Figure~\ref{Fig3D}, where $f$ 
is the characteristic function of   $[0.15,0.3]\cup [0.5,0.75]$. 
We study the resilience of the super-resolution demonstrated in these examples with respect to various levels of noise.
 In all our examples,  one can see that selecting the regularisation parameter $m$  according to the  residual minimisation  principle is adequate. In particular, 
 we observe that the PSWF-Radon reconstruction   $\tilde{f}_m$ with such $m$ performs at least as well as $\naive$, even for high levels of noise, while $\tilde{f}_m$   significantly outperforms $\naive$ for low levels of noise.
 In contrast to Figure~\ref{Fig2D} and Figure~\ref{Fig3D}, we present our numerical results as 1D plots instead of 2D images, as this format is easier for comparisons.

\subsubsection{Order $\nu=0$}\label{S:two-steps0}

\begin{figure}[H]
\centering
 \begin{subfigure}{0.31\textwidth}
 \includegraphics[height=3.3cm]{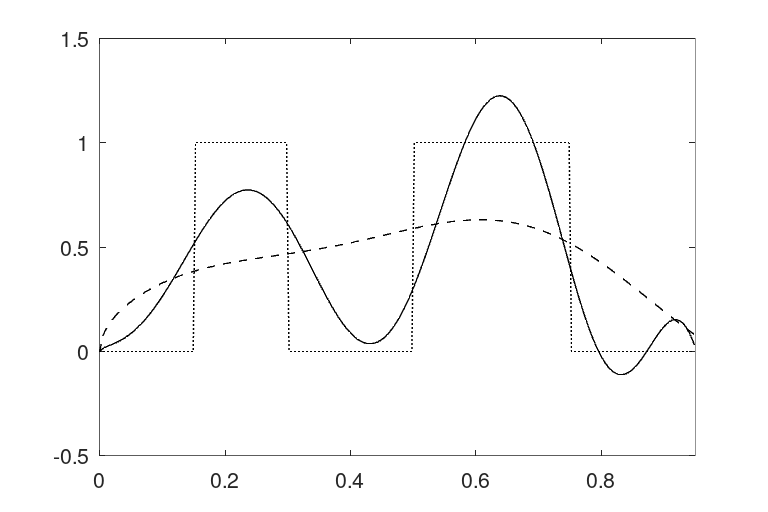}
 \caption{
 }
\end{subfigure}
 \begin{subfigure}{0.31\textwidth}
 \includegraphics[height=3.3cm]{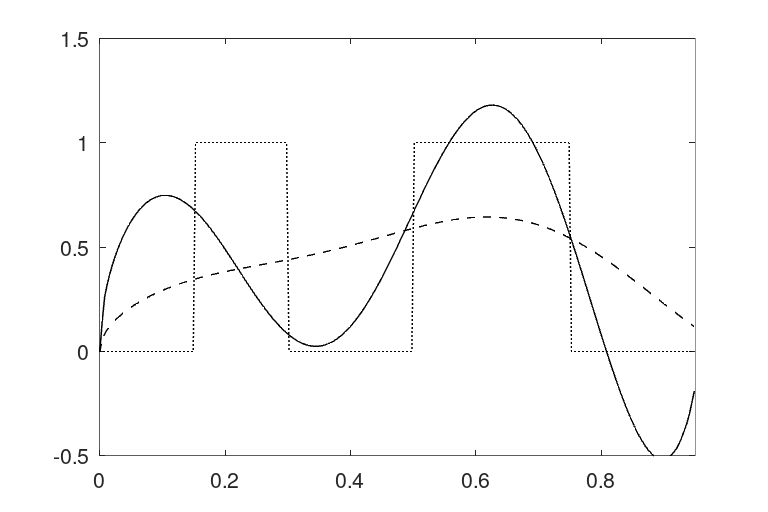}
 \caption{
 }
\end{subfigure}
 \begin{subfigure}{0.31\textwidth}
 \includegraphics[height=3.3cm]{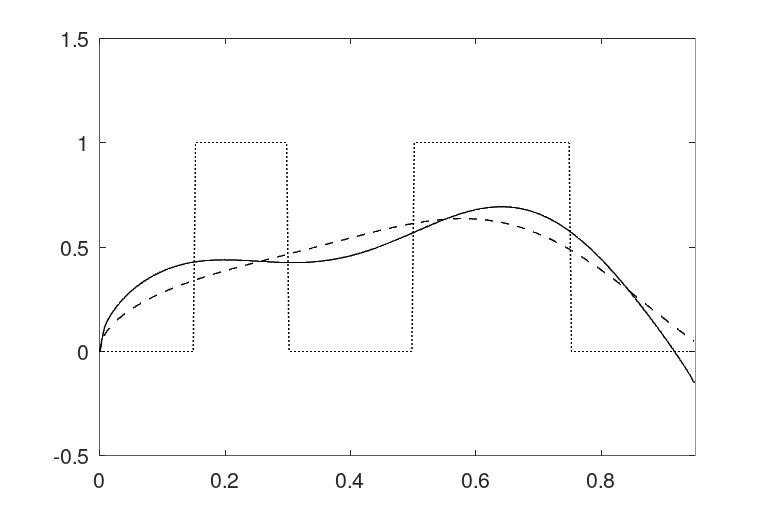}
 \caption{
 }
\end{subfigure}
\captionsetup{singlelinecheck=off}
 \caption[foo bar]{ \label{Fig:two-step-config0} 
   The PSWF-Radon reconstruction $\tilde{f}$    in comparison with   $\naive$ and the preimage $f$  for the order $\nu=0$:  (a) noiseless  
  (b) 20\% noise (c) 35\% noise.
 }
\end{figure}

\begin{figure}[H]
\centering
 \begin{subfigure}{0.31\textwidth}
 \includegraphics[height=3.3cm]{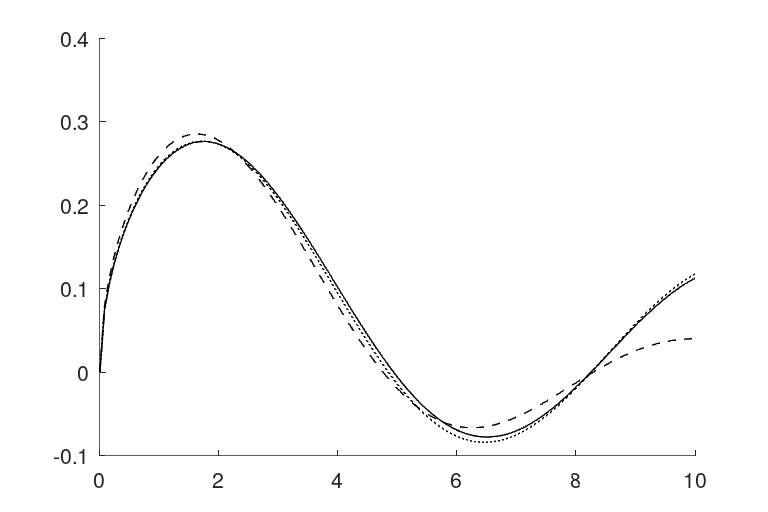}
 \caption{
 }
\end{subfigure}
 \begin{subfigure}{0.31\textwidth}
 \includegraphics[height=3.3cm]{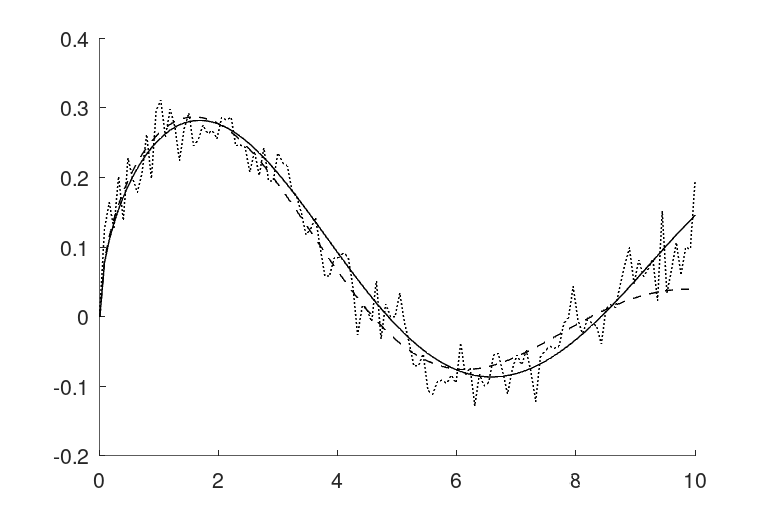}
 \caption{
 }
\end{subfigure}
 \begin{subfigure}{0.31\textwidth}
 \includegraphics[height=3.3cm]{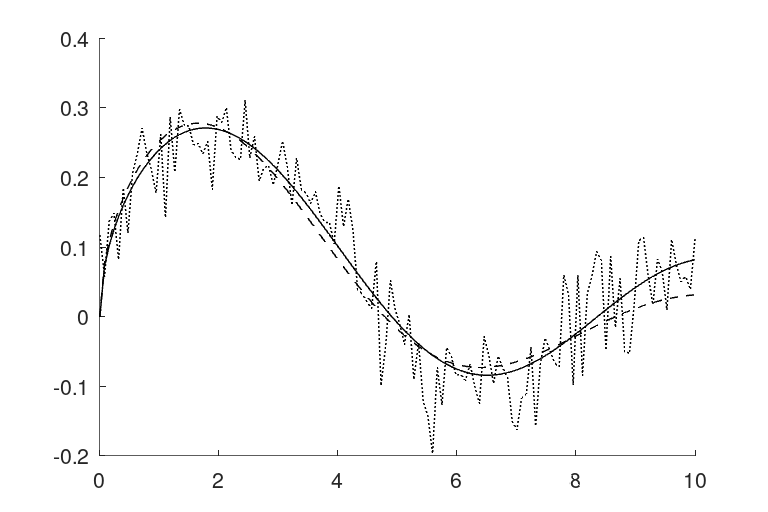}
 \caption{
 }
\end{subfigure}
\captionsetup{singlelinecheck=off}
 \caption[foo bar]{ 
 \label{Fig:two-step-data0}
 The Hankel transforms 
    $\calH_\nu[\tilde{f}]$      
    and $\calH_\nu[\naive]$    
    in comparison with  the Hankel data $h$  for the order $\nu =0$:  (a) noiseless,  
  (b) 20\% noise, (c) 35\% noise.
 }
\end{figure}

\begin{figure}[H]
\centering
 \begin{subfigure}{0.31\textwidth}
 \includegraphics[height=3.3cm]{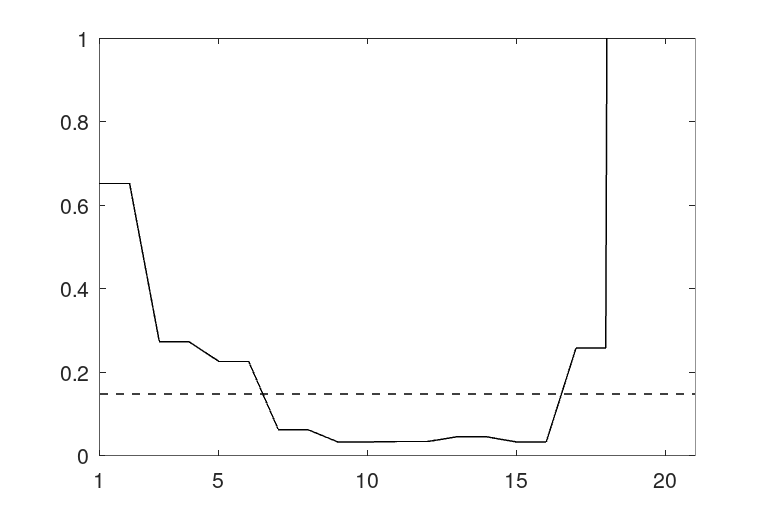}
 \caption{
 }
\end{subfigure}
 \begin{subfigure}{0.31\textwidth}
 \includegraphics[height=3.3cm]{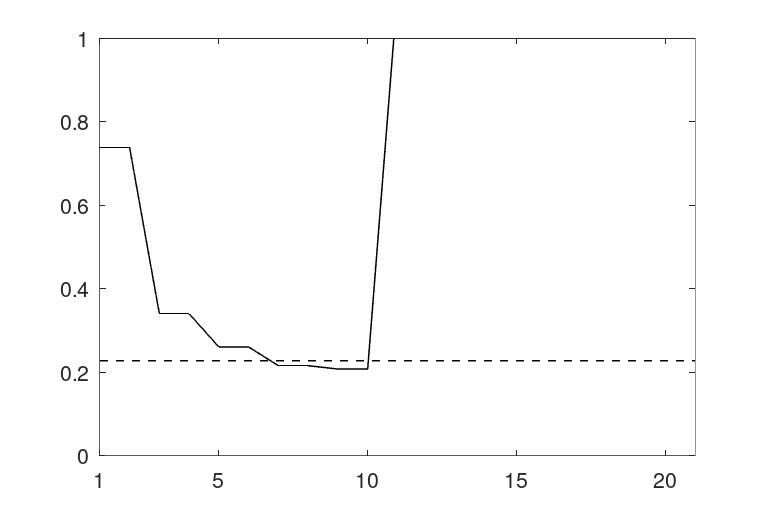}
 \caption{
 }
\end{subfigure}
 \begin{subfigure}{0.31\textwidth}
 \includegraphics[height=3.3cm]{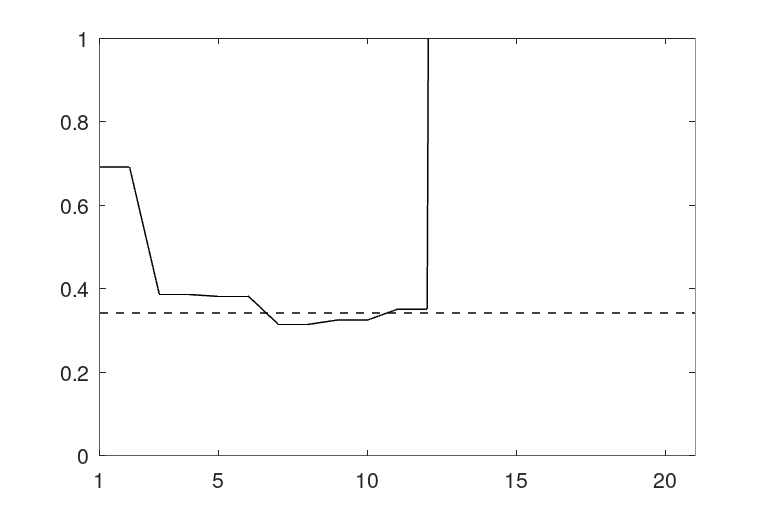}
 \caption{
 }
\end{subfigure}
\captionsetup{singlelinecheck=off}
 \caption[foo bar]{ 
  \label{Fig:two-step-residuals0}
 The dependence of the relative residual $\Err( \tilde{f}_m, h)$   on $m$ in comparison with  $\Err(\naive,h)$  for the order $\nu=0$: 
(a) noiseless,  
  (b) 20\% noise, (c) 35\% noise.
 }
\end{figure}

Here, we consider the PSWF-Radon reconstruction $\tilde{f}$ for $\nu=0$ from its Hankel data 
given with three levels of white Gaussian noise: 0\%, 20\%, and 35\%.  
Figure~\ref{Fig:two-step-config0}(a) is equivalent to 
Figure~\ref{Fig2D} and illustrates the reconstruction $\tilde{f}$ from noiseless data.
Figure~\ref{Fig:two-step-config0}(b,c) illustrate the reconstructions $\tilde{f}$ from noisy data:  
one can see  that $\tilde{f}$
maintains super-resolution  even for $20\%$ noise (even though the location of the first step starts to shift), whilst   
$\tilde{f}$ resembles $\naive$ for 35\% noise.  
The corresponding plots in the data space are given in  Figure  
\ref{Fig:two-step-data0}. 
  Figure \ref{Fig:two-step-residuals0} illustrates the evolution of the residual plots of $\Err( \tilde{f}_m, h)$ defined in \eqref{def:Err}: as noise increases, the residual minimisation   leads to smaller  $m$ and eventually the PSWF-Radon reconstruction $\tilde{f}$ loses its capability to super-resolution  and becomes comparable to  $\naive$.

\subsubsection{Order $\nu=0.5$}\label{S:two-steps05}
\begin{figure}[H]
\centering
 \begin{subfigure}{0.31\textwidth}
 \includegraphics[height=3.3cm]{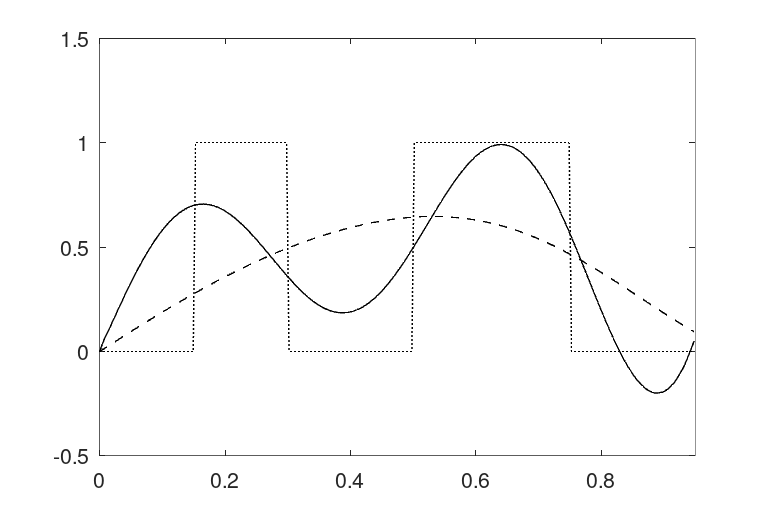}
 \caption{
 }
\end{subfigure}
 \begin{subfigure}{0.31\textwidth}
 \includegraphics[height=3.3cm]{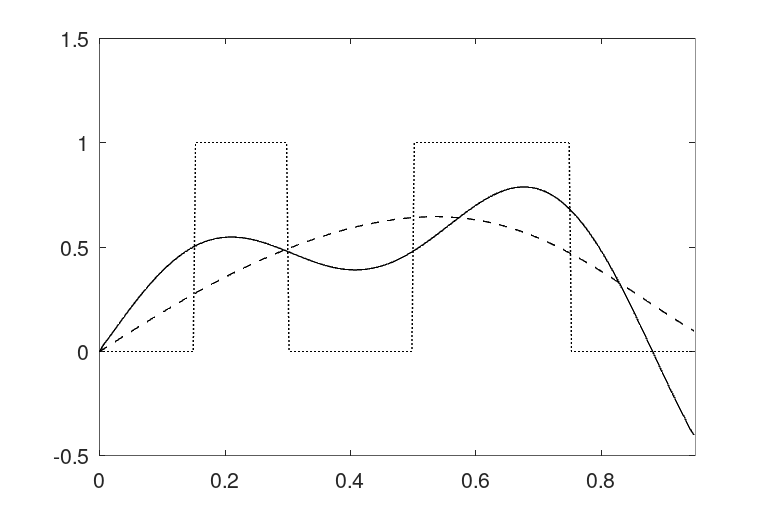}
 \caption{
 }
\end{subfigure}
 \begin{subfigure}{0.31\textwidth}
 \includegraphics[height=3.3cm]{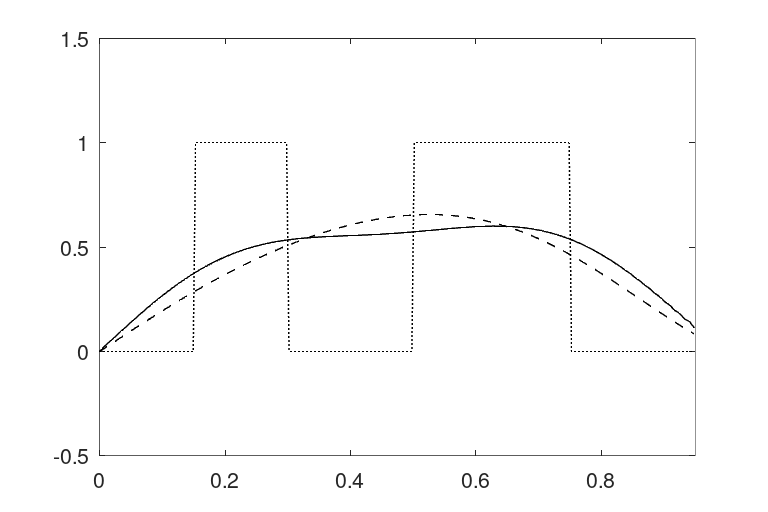}
 \caption{
 }
\end{subfigure}
\captionsetup{singlelinecheck=off}
 \caption[foo bar]{ 
 \label{Fig:two-step-config05}
    The PSWF-Radon reconstruction $\tilde{f}$    in comparison with   $\naive$   and the preimage $f$ for the order $\nu=0.5$:  (a) noiseless,  
  (b) 5\% noise, (c) 10\% noise.
 }
\end{figure}

\begin{figure}[H]
\centering
 \begin{subfigure}{0.31\textwidth}
 \includegraphics[height=3.3cm]{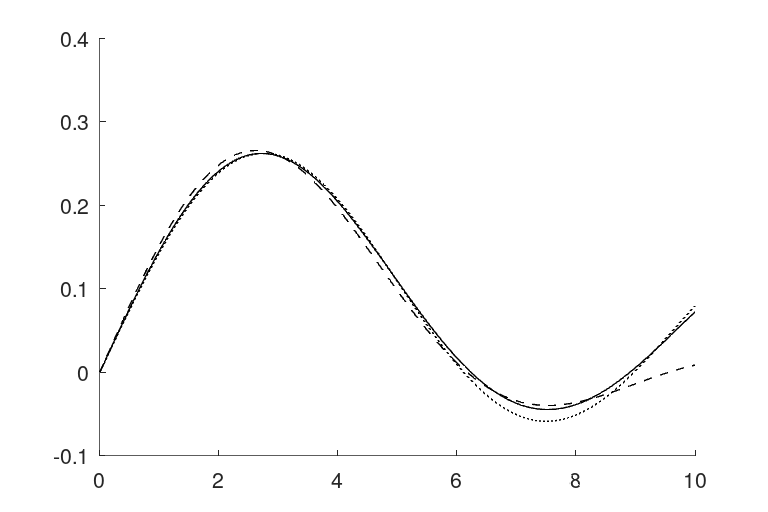}
 \caption{
 }
\end{subfigure}
 \begin{subfigure}{0.31\textwidth}
 \includegraphics[height=3.3cm]{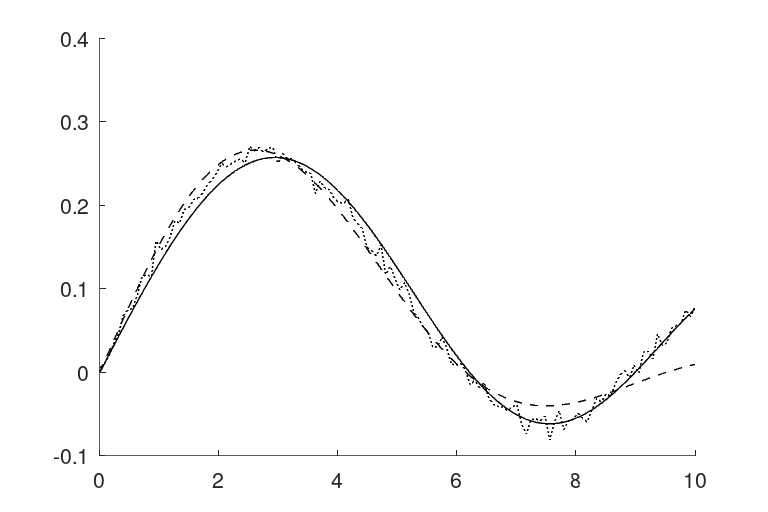}
 \caption{
 }
\end{subfigure}
 \begin{subfigure}{0.31\textwidth}
 \includegraphics[height=3.3cm]{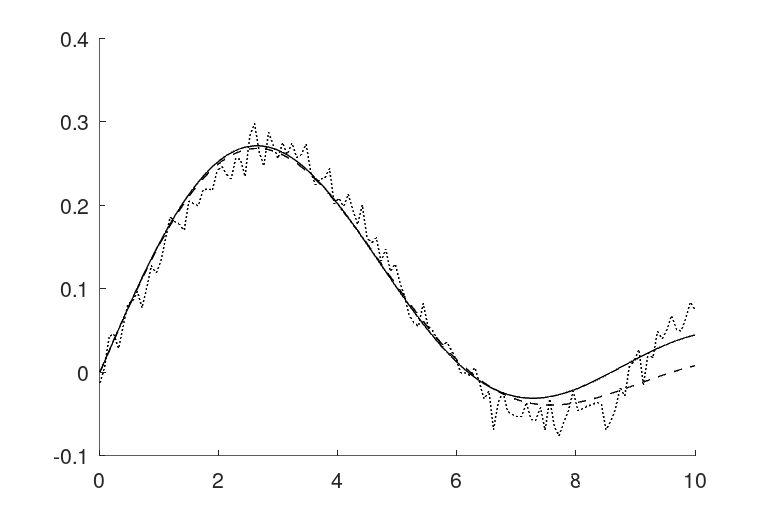}
 \caption{
 }
\end{subfigure}
\captionsetup{singlelinecheck=off}
 \caption[foo bar]{ 
 \label{Fig:two-step-data05}
The Hankel transforms 
    $\calH_\nu[\tilde{f}]$      
    and $\calH_\nu[\naive]$    
    in comparison with  the Hankel data $h$ for the order $\nu=0.5$:  (a) noiseless,  
  (b) 5\% noise, (c) 10\% noise.
 }
\end{figure}

Here, we consider the PSWF-Radon reconstruction $\tilde{f}$ for $\nu=0.5$ from its  Hankel data 
given with three levels of white Gaussian noise: 0\%, 5\%, and 10\%.  
Figure~\ref{Fig:two-step-config05}(a) is equivalent to 
Figure~\ref{Fig3D} and illustrates the reconstruction $\tilde{f}$ from noiseless data.
Figure~\ref{Fig:two-step-config05}(b,c) illustrate the reconstructions $\tilde{f}$ from noisy data:  
one can see  that $\tilde{f}$
maintains super-resolution   
for $5\%$ noise 
(which is lower than for the case of $\nu=0$), whilst   
$\tilde{f}$ resembles $\naive$ already for 10\% noise.  
The corresponding plots in the data space are given in  Figure  
\ref{Fig:two-step-data05}. 
 Figure~\ref{Fig:two-step-config05}
 illustrates the evolution of the residual plots of $\Err( \tilde{f}_m, h)$ 
  similarly to Figure~\ref{Fig:two-step-residuals0}.

\begin{figure}[H]
\centering
 \begin{subfigure}{0.31\textwidth}
 \includegraphics[height=3.3cm]{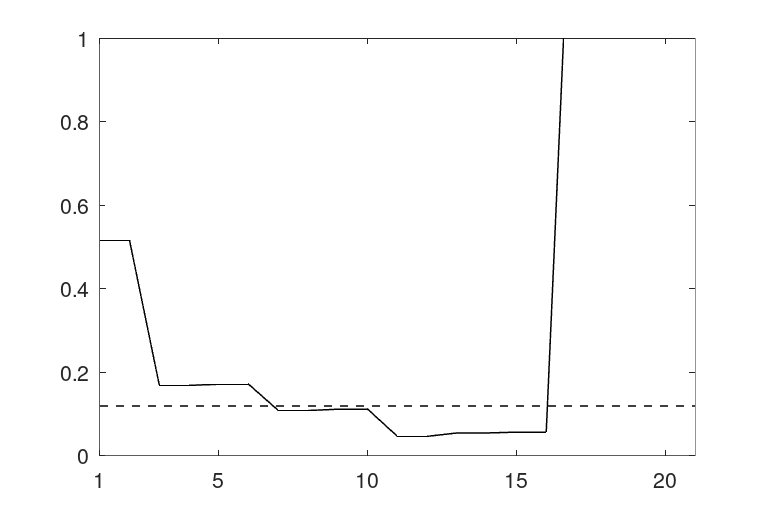}
 \caption{
 }
\end{subfigure}
 \begin{subfigure}{0.31\textwidth}
 \includegraphics[height=3.3cm]{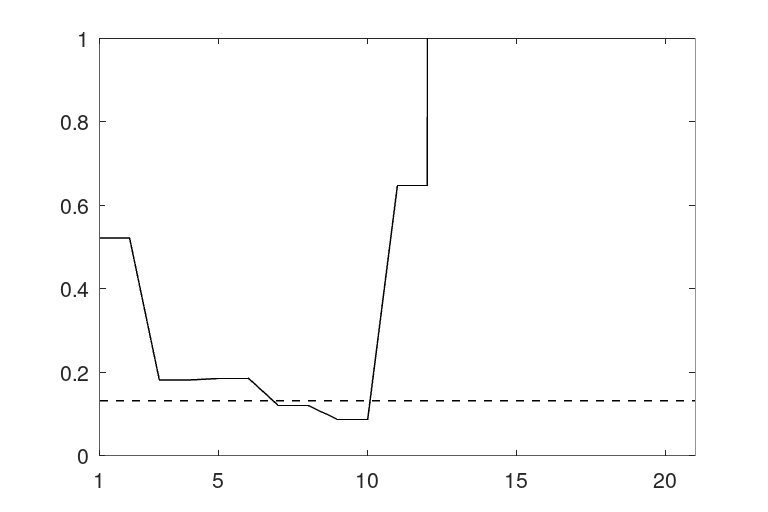}
 \caption{
 }
\end{subfigure}
 \begin{subfigure}{0.31\textwidth}
 \includegraphics[height=3.3cm]{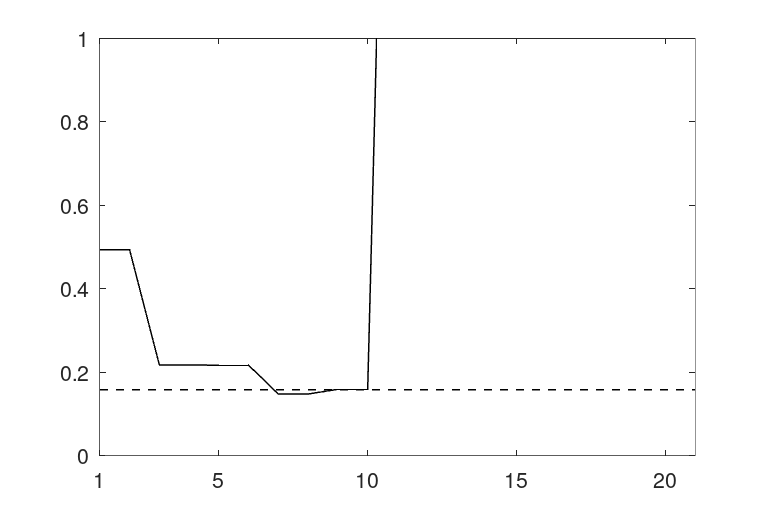}
 \caption{
 }
\end{subfigure}
\captionsetup{singlelinecheck=off}
 \caption[foo bar]{ 
 \label{Fig:two-step-residuals05}
 The  dependence of the relative residual $\Err( \tilde{f}_m, h)$   on $m$ in comparison with  $\Err(\naive,h)$  for the order $\nu=0.5$: 
(a) noiseless,  
  (b) 5\% noise, (c) 10\% noise.
 }
\end{figure}


 \subsection{Harmonics}\label{S:harmonics}
 Here, we  consider 
$f =  \sin (\omega s)$ truncated to the interval 
$s\in [0,\sigma]$, where $\sigma=1$, from its Hankel data limited to $[0,r]$, where $r=10$. This example is inherently difficult for Problem \ref{P1}  
as the Hankel data on $[0,r]$  rapidly decays shortly after  $\omega$ exceeds $r$ (that is, it becomes  an insignificant part of the  full Hankel transform).  In particular, $\naive$ fails to reconstruct $f$ adequately for such $\omega$.   
We examine the extent to which
the  PSWF-Radon reconstruction $\tilde f$ can overcome this limit for  the orders $\nu = 0$  and $\nu= 0.5$ at various  levels of noise.

  \subsubsection{Order $\nu= 0$}

\begin{figure}[H]
\centering
 \begin{subfigure}{0.31\textwidth}
 \includegraphics[height=3.3cm]{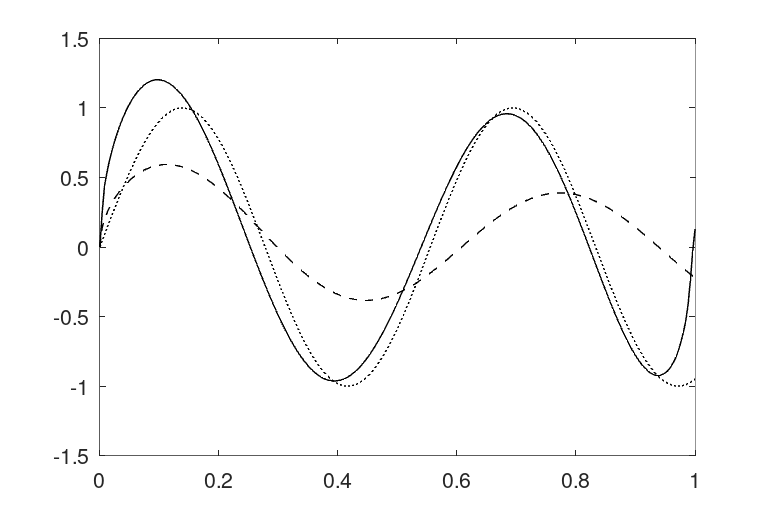}
 \caption{
 }
\end{subfigure}
 \begin{subfigure}{0.31\textwidth}
 \includegraphics[height=3.3cm]{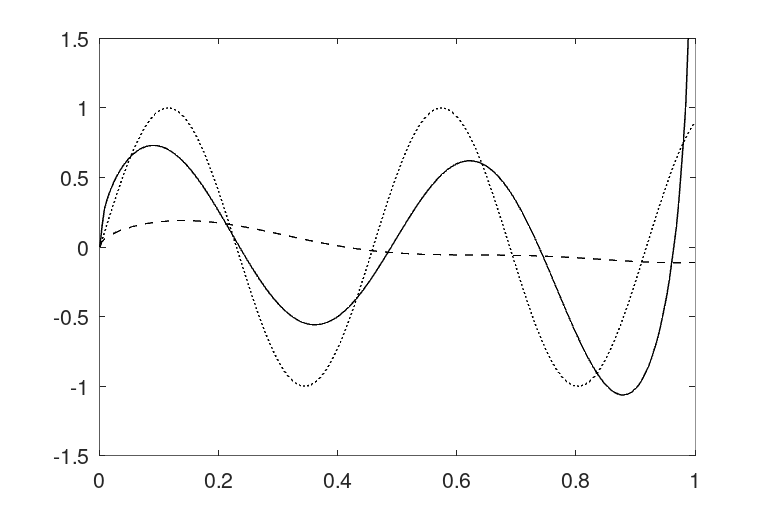}
 \caption{
 }
\end{subfigure}
 \begin{subfigure}{0.31\textwidth}
 \includegraphics[height=3.3cm]{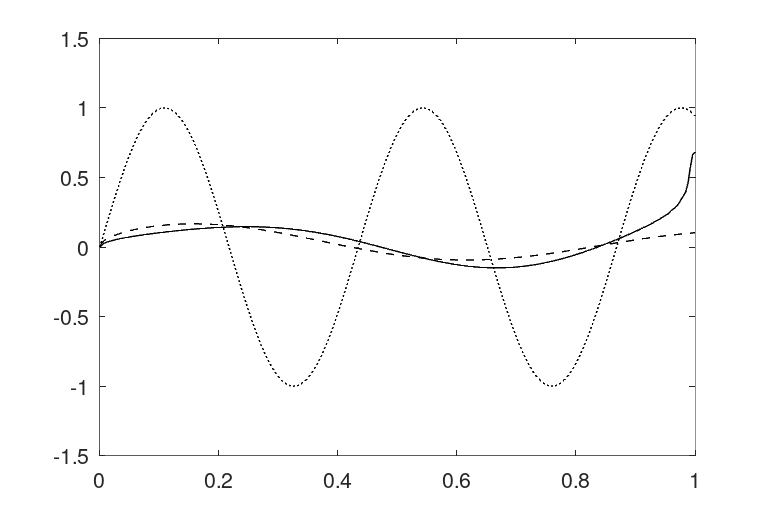}
 \caption{
 }
\end{subfigure}
\captionsetup{singlelinecheck=off}
 \caption[foo bar]{ \label{F:harmonics0_20}
   The PSWF-Radon reconstruction $\tilde{f}$     in comparison with   $\naive$  and the preimage $f=\sin(\omega x)$  
   for the order $\nu=0$ and the Hankel data with 20\% noise:  (a)~$\omega= 11.32$,  
  (b)~$\omega=13.68$,  (c)~$\omega=14.47$.
 }
\end{figure}

Figure \ref{F:harmonics0_20} shows the evolution with respect to $\omega$ of the reconstructions $\tilde{f}$ and $\naive$  for the case of   20\% white Gaussian noise and $\nu=0$. Figure \ref{F:harmonics0_20}(a) 
illustrates that  $\naive$
has already deteriorated significantly from $f$ at frequency $\omega =11.32$, whereas $\tilde{f}$ is fairly accurate.
 Figure \ref{F:harmonics0_20}(b) 
shows that  $\tilde{f}$ starts to deteriorate from $f$ around $\omega =13.68$, whereas $\naive$ essentially reduces to $0$.  Figure \ref{F:harmonics0_20}(c)  illustrates the behaviour of these reconstructions at   higher frequencies. Thus, for this level of noise, the  PSWF-Radon approach  extends the range of frequencies by around $20\%$ in comparison with $\naive$.

\begin{figure}[H]
\centering
 \begin{subfigure}{0.31\textwidth}
 \includegraphics[height=3.3cm]{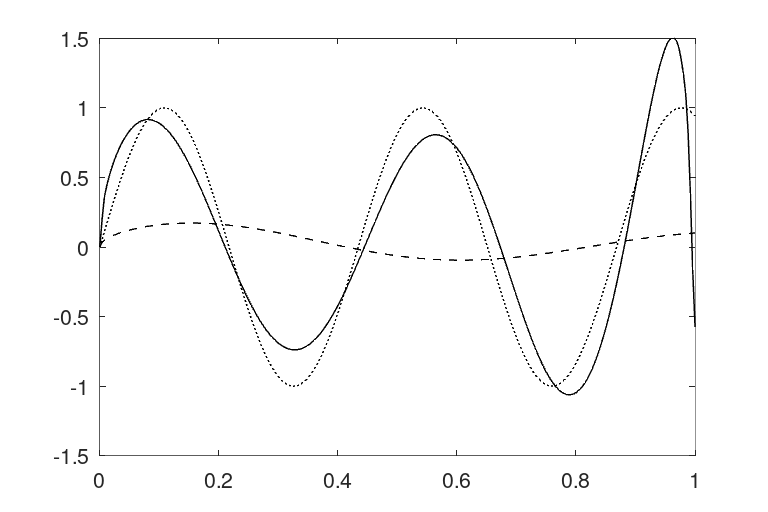}
 \caption{
 }
\end{subfigure}
 \begin{subfigure}{0.31\textwidth}
 \includegraphics[height=3.3cm]{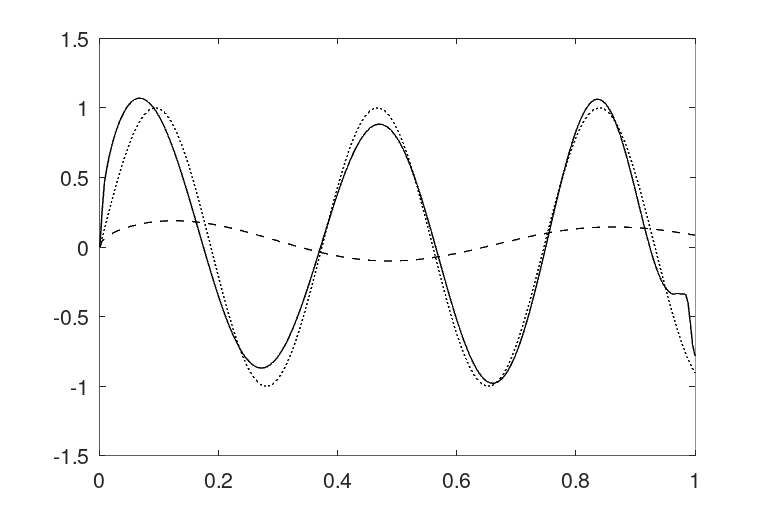}
 \caption{
 }
\end{subfigure}
 \begin{subfigure}{0.31\textwidth}
 \includegraphics[height=3.3cm]{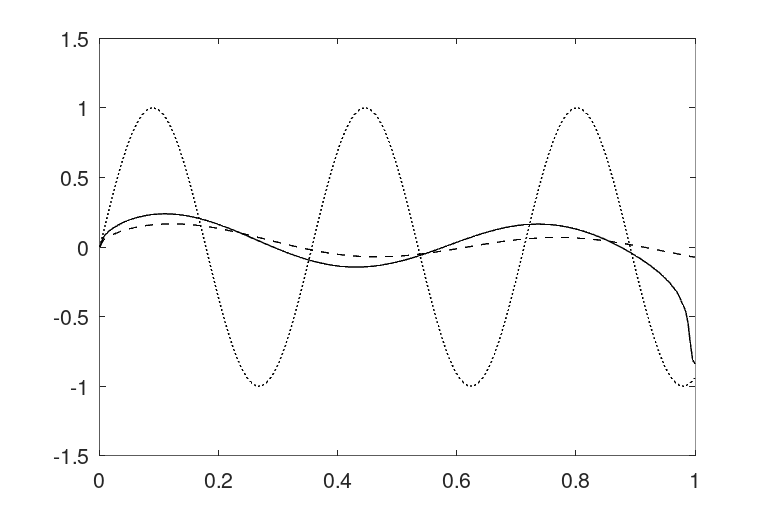}
 \caption{
 }
\end{subfigure}
\captionsetup{singlelinecheck=off}
 \caption[foo bar]{ \label{F:harmonics0}
   The PSWF-Radon reconstruction $\tilde{f}$   in comparison with   $\naive$   and the preimage $f=\sin(\omega x)$   for the order  $\nu=0$ and noiseless Hankel data:  (a)~$\omega=14.47$,  
  (b)~$\omega=16.84$,  (c)~$\omega=17.63$.
 }
\end{figure}

In turn, for the noiseless case, we  succeeded to extend  the range of frequencies by around $50\%$ in comparison with $\naive$.
Figure \ref{F:harmonics0} shows the evolution with respect to $\omega$ of the PSWF-Radon reconstruction $\tilde{f}$  for the noiseless case beyond the frequencies considered in  Figure~\ref{F:harmonics0_20}. Note that  $\naive$ remains essentially   $0$ for such high frequencies.    Figure~\ref{F:harmonics0}(a,b) shows that  $\tilde{f}$ is fairy accurate for frequencies up to around 17.  Figure~\ref{F:harmonics0}(c) shows that, in practice, the PSWF-Radon approach is limited by the precision of numerical computations, despite its theoretical capability to handle arbitrary frequencies.   

\subsubsection{Order $\nu = 0.5$}

\begin{figure}[H]
\centering
 \begin{subfigure}{0.31\textwidth}
 \includegraphics[height=3.3cm]{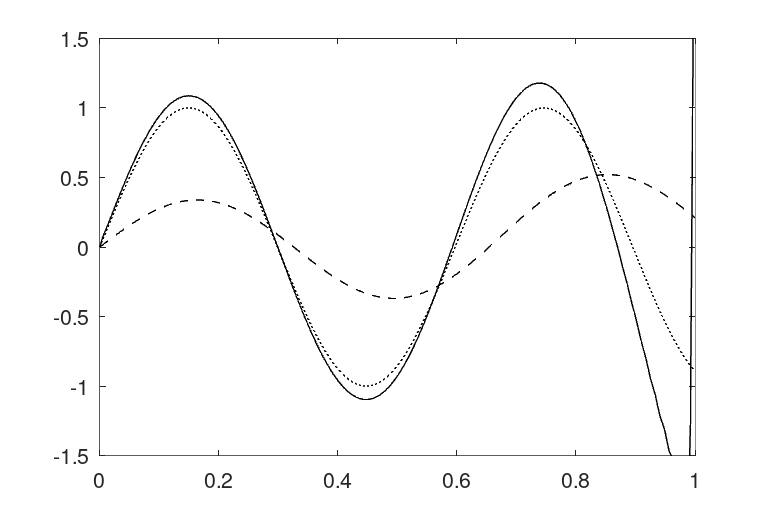}
 \caption{
 }
\end{subfigure}
 \begin{subfigure}{0.31\textwidth}
 \includegraphics[height=3.3cm]{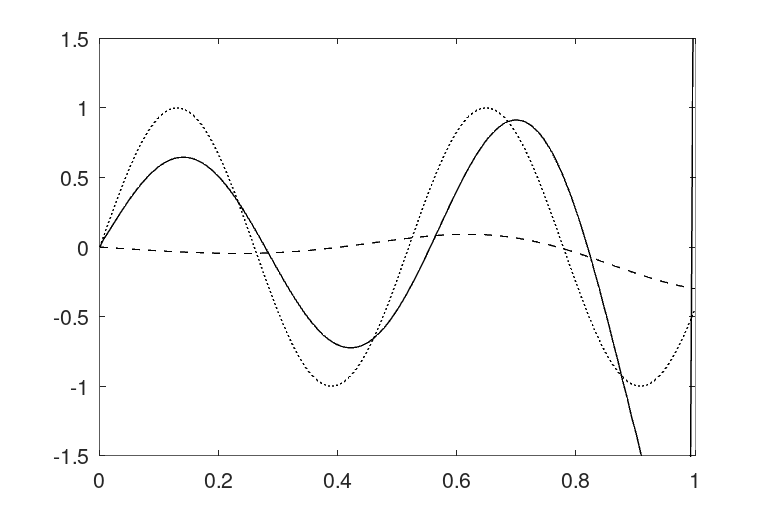}
 \caption{
 }
\end{subfigure}
 \begin{subfigure}{0.31\textwidth}
 \includegraphics[height=3.3cm]{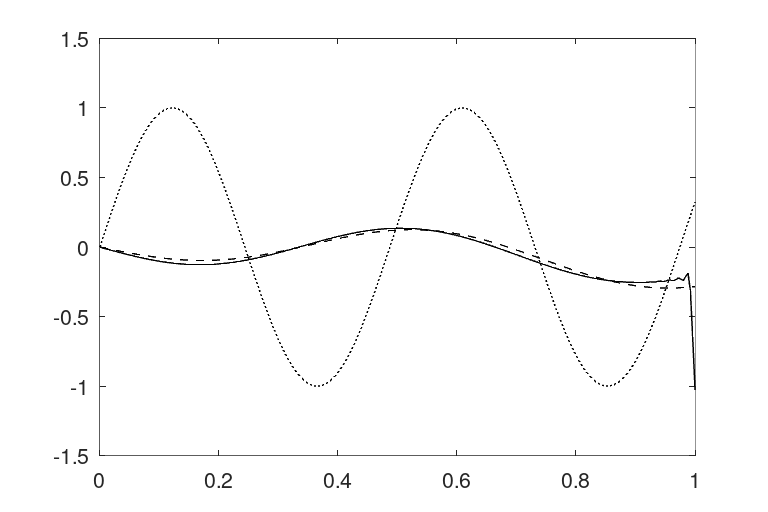}
 \caption{
 }
\end{subfigure}
\captionsetup{singlelinecheck=off}
 \caption[foo bar]{ \label{F:harmonics05_5}
   The  PSWF-Radon reconstruction $\tilde{f}$     in comparison with   $\naive$   and the preimage $f=\sin(\omega x)$  
   for the order $\nu=0.5$ and  the Hankel data with 5\% noise:  (a)~$\omega= 10.53$,  
  (b)~$\omega=11.32$,  (c)~$\omega=12.11$.
 }
\end{figure}

Figure \ref{F:harmonics05_5} shows the evolution with respect to $\omega$ of the reconstructions $\tilde{f}$ and $\naive$  for the case of   5\% white Gaussian noise and $\nu=0.5$. Figure \ref{F:harmonics05_5}(a) 
illustrates that  $\naive$
has already deteriorated significantly from $f$ at frequency $\omega =10.53$, whereas $\tilde{f}$ is fairly accurate.
 Figure \ref{F:harmonics05_5}(b) 
shows that  $\tilde{f}$ starts to deteriorate from $f$ around $\omega =11.32$, whereas $\naive$ essentially reduces to $0$.  Figure \ref{F:harmonics05_5}(c)  illustrates the behaviour of these reconstructions at   higher frequencies. Thus, for this level of noise, the  PSWF-Radon approach  extends the range of frequencies by around $7.5\%$ in comparison with $\naive$.

\begin{figure}[H]
\centering
 \begin{subfigure}{0.31\textwidth}
 \includegraphics[height=3.3cm]{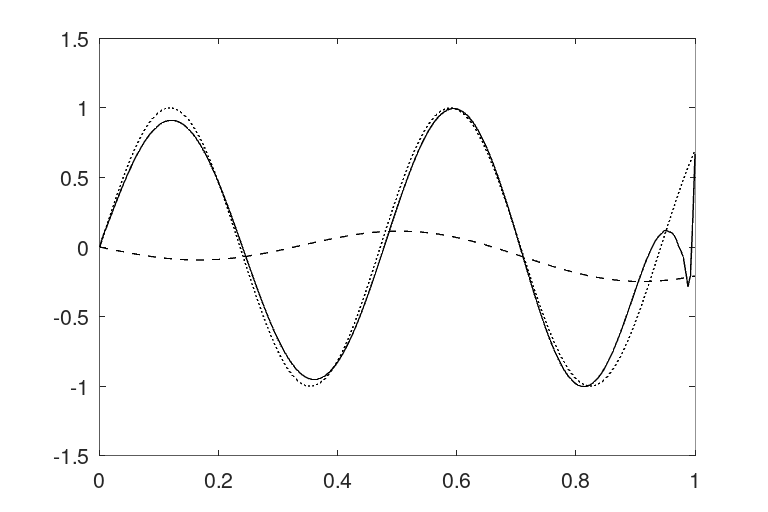}
 \caption{
 }
\end{subfigure}
 \begin{subfigure}{0.31\textwidth}
 \includegraphics[height=3.3cm]{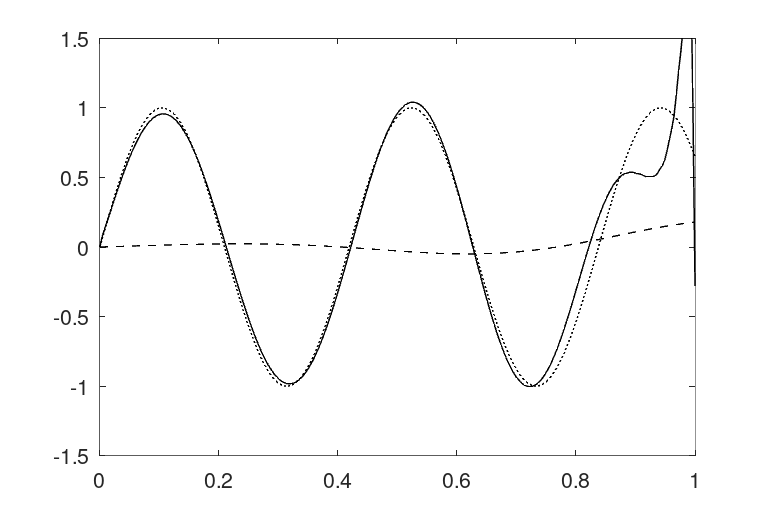}
 \caption{
 }
\end{subfigure}
 \begin{subfigure}{0.31\textwidth}
 \includegraphics[height=3.3cm]{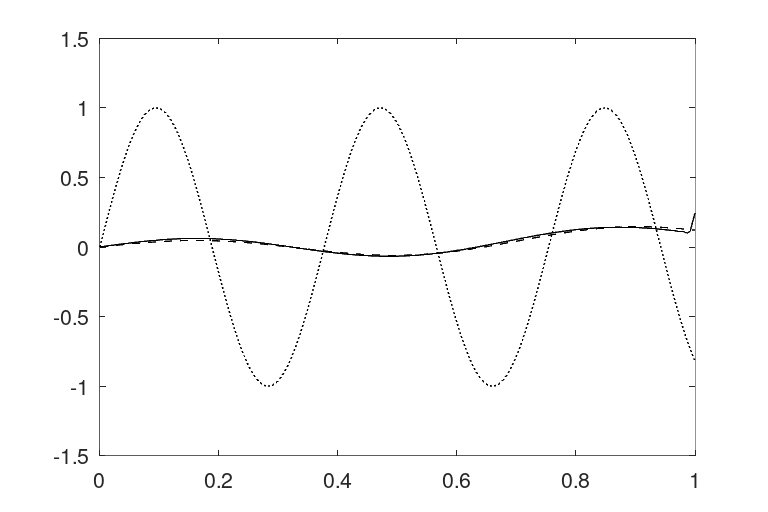}
 \caption{
 }
\end{subfigure}
\captionsetup{singlelinecheck=off}
 \caption[foo bar]{ \label{F:harmonics05}
   The PSWF-Radon reconstruction $\tilde{f}$     in comparison with   $\naive$   and the preimage $f=\sin(\omega x)$   for the order  $\nu=0.5$ and  noiseless Hankel data:  (a)~$\omega=13.33$,  
  (b)~$\omega=15$,  (c)~$\omega=16.67$.
 }
\end{figure}

In turn, for the noiseless case, we  succeeded to extend  the range of frequencies by around $40\%$ in comparison with $\naive$.
Figure \ref{F:harmonics05} shows the evolution with respect to $\omega$ of the PSWF-Radon reconstruction $\tilde{f}$  for the noiseless case beyond the frequencies considered in  Figure~\ref{F:harmonics05_5}. Note that  $\naive$ remains essentially   $0$ for such high frequencies.    Figure~\ref{F:harmonics05}(a,b) shows that  $\tilde{f}$ is fairy accurate for frequencies up to around 15.  Figure~\ref{F:harmonics05}(c)    illustrates the same conclusion as for Figure~\ref{F:harmonics0}(c).

 From the figures presented in Section \ref{S:harmonics}, one can  see that, in the case of noiseless data,    the PSWF-Radon reconstructions $\tilde{f}$  
 for  $\nu=0$ and $\nu=0.5$ have similar 
  super-resolution capability.
However,  the performance of $\tilde{f}$ for $\nu=0$ is notably superior in the presence of noise.



\subsection{Implementation details}
We follow Steps 1-3 formulated at the beginning of Section~\ref{S:numerics}.

In our examples for Problem \ref{P1},  the Hankel data $h$ is given on the uniform grid with  256 points on interval $[0,r]$, with $r=10$.  We 
reconstruct unknown 
function $f$ on uniform grid with 256 points on $[0,\sigma]$, where $\sigma=1$.    To generate data and also  compute the residuals in the data space, we evaluate  integral of \eqref{def:Hankel} using  the trapezoidal rule.

Our numerical implementation of  operator   $\calF_{m,c}^{-1}$  defined in \eqref{def:Fnc} relies on the Spheroidal Library of \cite{AGD2014}, which is available at 
\texttt{https://github.com/radelman/scattering}.  Unfortunately, eigenvalues of the PSWFs are not available in this library, so we computed them using the approach of \cite[Section 4.1]{HRY2001}.   It is important to note that the computations of the eigenvalues and the integrals in  \eqref{def:Fnc} require oversampling due to high oscillations of PSWFs.
Using linear interpolation for the Hankel data, we get samples on the uniform grid with 1024 points.

In dimension $d=2$, which corresponds to integer $\nu$, we use the Fourier-based implementation of filtered back projection (FBP) algorithm for the Radon transform inversion
from  NFFT3 library; see  survey  \cite{Keiner2009using}.
 In dimension $d=3$, which corresponds to   half-integer $\nu$,
 we also use this library  for reconstructions from the integrals over planes. Note that this version of the FBP algorithm is significantly less common so we relied on the numerical implementations from   \cite{Goncharov2021}.

Note that at Step~2, for fixed $m$, only one inversion by $\calF^{-1}_{m,c}$ is sufficient in view of separation of radial and angle variables $x$ and $\theta$, respectively. Therefore, the overall memory and time complexity of the algorithm is dominated by the Radon inversion. For fixed $\nu\in \mathbb{N}$, reconstruction was performed on a personal laptop (AMD Ryzen 7 7840U) using 16 threads and required less than 1 second (for grid sizes up to 512) with negligible memory usage. For $\nu\in \frac{1}{2} + \mathbb{N}$, reconstruction required nearly 2 minutes on workstation (AMD 3995W) with 64 threads and 60GB of RAM. Such difference of required resources is due to dimensionality of the corresponding Radon inversion -- 2D and 3D for integer and half-integer $\nu$, respectively.

In our experiments, we also observed the Gibbs phenomenon at the boundary point $x=\sigma = 1$ for reconstruction $\tilde{f}$, see, for example, Figure~\ref{F:harmonics05_5}(c)
and Figure~\ref{F:harmonics05}. 
We explain this phenomenon by a combination of the following factors: 
\begin{itemize}
    \item harmonic phantoms in Figures~\ref{F:harmonics05_5},~\ref{F:harmonics05} do not vanish at the boundary point $s=1$;
    \item function $w_{r,\sigma}$ in Theorem~\ref{T:Fourier} reconstructed for fixed $m$ using $\mathcal{F}_{m,c}^{-1}$ does not necessarily belong to the range of the Radon transform on functions supported in $B_1$;
    \item our implementation of Radon inversion strongly relies on the assumption of compactly supported input (we use \textit{zero padding} of the input to increase precision of the corresponding Fourier integration).
\end{itemize}
In particular, recall that the Radon transform $w_{r,\sigma}(x,\theta)$ of a function supported in the ball $B_1$ vanishes as $|x| \rightarrow 1$ at a specific rate (see \cite[Theorem 1.6]{Naterrer2001}). 
To decrease the ripple at the boundary, one      may   increase $\sigma$ so that considered $f=f(s)$ in Problem~\ref{P1} vanishes near the boundary point $s=\sigma$ or/and apply a smooth cut-off for $w_{r,\sigma}(x,\theta)$ near the boundary $|x|=1$. 
However, we did not pursue this task as auxiliary regularisation might affect the investigated stability.

The source code of our algorithms based on the PSWF-Radon approach  generating the presented data can be found   at \texttt{https://github.com/fedor-goncharov/pswf-radon}.



\section{Concluding remarks}\label{S:concluding}

 Note that   Problem \ref{P1} is equivalent to the inversion of the band-limited Hankel transform $\calH_{\nu,c}$  with $c= r \sigma$ defined by  
\begin{equation*} 
    \calH_{\nu,c}[f](x) 
:= \int_0^{1} f(y) J_{\nu}(cxy) \sqrt{cxy}\, d y, \qquad x \in [0,1].
\end{equation*}
The operator $\calH_{\nu,c}$ is well-studied in the literature and its properties are similar in many respects to the band-limited Fourier operator $\calF_c$;  see, for example,  \cite{KM2009,Slepian1964}. 
One can approach Problem \ref{P1} 
using the singular value decomposition for the operator $\calH_{\nu,c}$. 
 Compared to this approach, the main advantages of the present work based on  the PSWF-Radon approach of \cite{INS2022,INS_Proceedings,INnotePSWF} are summarised below.   
 \begin{itemize}
 \item Theorem \ref{T:int} and Theorem \ref{T:half} involve classical functions only, namely, the Chebyshev and Legendre polynomials  and the prolate spheroidal wave functions of \cite{Niven1880, SP1961}.

 \item
 In particular, our formulas  require the same set of eigenfunctions $(\psi_{j,c})_{j \in \Naturals}$ for different $\nu$, that is, classical PSWFs,  unlike the approach 
 based on the eigenfunctions of $\calH_{\nu,c}$. This is particularly convenient  for problems  involving simultaneous inversion of $\calH_{\nu,c}$ for many~$\nu$.



 \item  
 We proceed from efficient numerical implementation  for  band-limited Fourier inversion in multidimensions   developed in \cite{INS2022, INS_Proceedings} that achieves super-resolution even for noisy data. 
 We observe similar numerical properties   for Problem~\ref{P1}.

 \end{itemize}

Next, we remark that our work gives the first numerical illustration of the PSWF-Radon approach for 
band-limited Fourier inversion in three dimensions; see Figure \ref{Fig3D} and Section \ref{S:two-steps05}.  
However, further numerical investigations are  required 
for reconstructions without a-priori assumption of spherical (or spherical-type) symmetry. 
An important challenge for the PSWF-Radon approach is that the observed resilience  of our reconstructions to the noise is significantly lower in the case of $d=3$ compared to $d=2$. 
Further natural challenges also include reduction of computational costs for PSWF-Radon reconstructions from band-limited Hankel data. 

 Furthermore, various a-priori assumptions on unknown function $f$ 
 in  Problem \ref{P1}  can enhance super-resolution properties. Possible approaches to this (and  similar improvements for band-limited Fourier inversion for $d=2$ and $d=3$) include combinations of the PSWF-Radon approach \cite{INnotePSWF, INS2022, INS_Proceedings} with the methods of \cite{AM,Meng2023} based on  generalised PSWFs,  and, for further enhancement, 
variational and iterative methods including machine learning  techniques; see, for example, \cite{AMS2009,SHN2023, LZBNJ2023} and references therein.

\section{Acknowledgements}

The work was initiated in the framework of  the internship of  R.~Zaytsev at  the Centre de Math\'ematique Appliqu\'ees of  Ecole Polytechnique
under the supervision of R.G.~Novikov  in June-August 2023.

\end{document}